\documentclass[leqno]{article}
\usepackage{amsmath}
\usepackage{amsfonts}
\usepackage{amssymb}
\usepackage{graphicx}

\newtheorem{theorem}{Theorem}

\newtheorem{lemma}[theorem]{Lemma}

\newtheorem{proposition}[theorem]{Proposition}
\newtheorem{remark}[theorem]{Remark}
\newtheorem{definition}[theorem]{Definition}
\newenvironment{proof}[1][Proof]{\noindent \textit{#1.} }{\ \rule{0.5em}{0.5em}}

\begin{document}

\title{\textbf{Connections functorially attached }\\
\textbf{to almost complex product structures}}
\author{Fernando Etayo\footnote{Departamento de Matem\'{a}ticas, Estad\'{\i}stica y Computaci\'{o}n.
 Facultad de Ciencias.  Universidad de Cantabria.
 Avda. de los Castros, s/n, 39071 Santander, SPAIN.
 e-mail: etayof@unican.es},  Rafael Santamar\'{\i}a\footnote{Departamento de Matem\'{a}ticas.
 Escuela de Ingenier\'{\i}as Industrial e
 Inform\'{a}tica.
 Universidad de Le\'{o}n.
 Campus de Vegazana, 24071 Le\'{o}n, SPAIN. e-mail: demrss@unileon.es}}
\date{}
\maketitle
\begin{abstract}
Manifolds endowed with three foliations pairwise transversal are known as 3-webs. Equivalently, they can be
algebraically  defined as biparacomplex or complex product manifolds, \emph{i.e.}, manifolds endowed with three
tensor fields of type $(1,1)$, $F$, $P$ and $J=F \circ P$, where the two first are product and the third one is
complex, and they mutually anti-commute. In this case, it is well known that there exists a unique torsion-free
connection parallelizing the structure. In the present paper, we study connections attached to non-integrable
  almost biparacomplex manifolds.
\end{abstract}

\bigskip\bigskip


\bigskip


\section{Introduction}
Theory of webs was begun by Blaschke \cite{BB}  when he introduced
a web on a surface as three families of curves pairwise
transversal at every point. Similarly, a three-web on a
$2n$-dimensional manifold $M$ is given by three $n$-dimensional
foliations such that their leaves through any point are pairwise
transverse. An $\alpha $-structure on $M$ is a system of three
distributions $V_{i}$ each of dimension $n$ such that
$TM=V_{1}\oplus V_{2}=V_{1}\oplus V_{3}=V_{2}\oplus V_{3}$.
Cruceanu obtained \cite{Cr} an algebraic characterization of these
$\alpha$-structures, by means of the almost biparacomplex
structures (see also \cite{San}).

An \emph{almost biparacomplex structure} on a manifold $M$ is given by two tensor fields $F$ and $P$ of type
$(1,1)$ satisfying $F^{2}=P^{2}=Id$, $F\circ P+P\circ F=0$. An almost biparacomplex structure is said to be
\emph{biparacomplex} if the distributions $T_{F}^{\pm}(M)=\ker(F\mp Id)$, $T_{P}^{\pm}(M)=\ker(P\mp Id)$
associated to the eigenvalues $\pm1$ of $F$ and $P$ respectively, are involutive. This condition is equivalent
to the vanishing of the Nijenhuis tensors $N_{F}$ and $N_{P}$.

An almost biparacomplex structure defines an
\emph{$\alpha$-structure} on the manifold by simply setting
$V_{1}=T_{F}^{+}(M)$, $V_{2}=T_{F}^{-}(M)$ and
$V_{3}=T_{P}^{+}(M)$. We call it the $\alpha$-structure associated
to the almost biparacomplex structure. In particular, a
biparacomplex structure defines a web. Conversely, given an
$\alpha$-structure $(V_{1},V_{2},V_{3})$, there exists a unique
almost biparacomplex structure such that $T_{F}^{+}(M)=V_{1}$,
$T_{F}^{-}(M)=V_{2}$, $T_{P}^{+}(M)=V_{3}$,
$T_{P}^{-}(M)=F(V_{3})$.

The distributions $T_{F}^{\pm}(M)$, $T_{P}^{\pm}(M)$  are readily seen to be equidimensional (and then $M$ is
endowed with two almost paracomplex structures, thus being called an almost biparacomplex manifold). Hence the
dimension of $M$ is even, and the tensor field $J=F\circ P$ defines an almost complex structure on $M$, so that
$M$ must be orientable.

Then an almost biparacomplex manifold is endowed with three tensor
fields of type $(1,1)$, $F$, $P$ and $J=F\circ P$, where the two
first are almost product and the third one is almost complex, and
they mutually anti-commute, thus being also called an \emph{almost
complex product manifold}. Such a structure was studied by
Libermann \cite{L} fifty years ago. In this sense, almost
biparacomplex manifolds corresponds to the \emph{paraquaternionic
numbers}, which consist on the 4-dimensional real algebra
generated by $\{ 1,i,j,k\}$ with the paraquaternionic relations:
$i^{2}=j^{2}=-k^{2}=1, \, ij=-ji=k$.
\vspace{3mm}

On the other hand, in recent years 3-webs have been considered in
a different  framework. We want to point out the works of Andrada,
Barberis, Bla\v{z}i\'{c}, Dotti, Ivanov, Kamada, Ovando,  Tsanov,
Vukmirovi\'{c}, Zamkovoy \cite{ABDO, AD, AS, BV, IT, ITZ, IZ,
Kamada} and others. Roughly speaking, they call a {\em complex
product structure} (also called a {\em para-hypercomplex
structure} and a \emph{hyper-paracomplex structure}) on a manifold
a pair of a complex structure $J$ and a product structure $P$ with
$J\circ P =-P\circ  J$. Of course, $F=J\circ P$ is a product
structure. They only consider the integrable case, \emph{i.e.},
the case where the Nijenhuis tensors of $J$ and $P$ (and that of
$F$) vanish. They know that there exists a unique torsion-free
connection parallelizing both structures $J$ and $P$. They mainly
apply their results to the study of Lie algebras.

\vspace{3mm}

In the present paper we shall study the general case, not only the integrable one, obtaining significant
advances:

\begin{itemize}

\item In any case, there exists a {\em canonical connection}
associated to the almost biparacomplex structure, defined as the
unique connection  parallelizing the tensor fields $F,P,J$ and
satisfying $T(X^{+},Y^{-})=0$, where $X^{+}\in
T_{F}^{+}(M),Y^{-}\in T_{F}^{-}(M)$  (Theorem
\ref{teor:cruceanu}), where $T$ denotes the torsion tensor field
of the connection.

\item The canonical connection is torsion-free iff the structure is integrable (Theorem \ref{teor:fn-torsion}).

\item The $G$-structure defined by the almost biparacomplex
structure is integrable iff the canonical connection is locally
flat (Theorem \ref{teor:cplana}).

\item The canonical connection is a functorial connection (Theorem
\ref{canonicaesfuntorial}).

\item In any case, there exists another functorial connection,
which is called the {\em well-adapted connection} (Theorem
\ref{existewelladapted}).

\item The canonical and the well-adapted connection coincide if
the structure is integrable (Theorem \ref{coinciden}).
\end{itemize}

As one can see these results enlightsheds light upon the connection defined in such a manifold and open some
questions about the non-integrable case in the aforementioned works.

The organization of the paper is as follows: in section 2 we
introduce notations. In section 3 we obtain the main results about
the canonical connection of an almost biparacomplex manifold. Also
we study functorial connections attached to such a manifold, in
particular the so called well-adapted connection. Finally, in the
last section we show examples, open problems and the relationship
among the quoted recent papers on complex product and
para-hypercomplex structures and the present paper.

\section{Notations and preliminaries\label{notandprel}}

Manifolds are assumed to be of class $C^{\infty}$ and satisfying
the second axiom of countability. Differentiable maps between
manifolds are also assumed to be of class $C^{\infty}$. We denote
by $\otimes^{k}V$, $\wedge^{k}V$, $S^{k}V$ the $k$th tensor,
exterior and symmetric power of a vector bundle $V$ over a
manifold $M$. In particular, we apply this notation to the tangent
and cotangent bundles $T(M)$, $T^{\ast}(M)$, respectively. We set
$\mathfrak{X}(M)=\Gamma(M,T(M))$,
$\Omega^{k}(M)=\Gamma(M,\wedge^{k}T^{\ast}(M))$, $k\in\mathbb{N}$.

If $\varphi\colon M\rightarrow M^{\prime}$ is a diffeomorphism,
for every $X\in\mathfrak{X}(M)$ we denote by $\varphi\cdot
X\in\mathfrak{X}(M^{\prime})$ the vector field defined by
$(\varphi\cdot X)_{x^{\prime}}=\varphi_{\ast }
(X_{\varphi_{-1}(x^{\prime})})$, $\forall x^{\prime}\in
M^{\prime}$. More generally, $\varphi$ transforms a tensor field
of type $(p,q)$ on $M$ into a tensor field of the same type on
$M^{\prime}$ by imposing
\[
\varphi\cdot(\omega_{1}\otimes\cdots\otimes\omega_{p}\otimes X_{1}
\otimes\cdots\otimes
X_{q})=\varphi^{-1\ast}\omega_{1}\otimes\cdots
\otimes\varphi^{-1\ast}\omega_{p}\otimes\varphi\cdot
X_{1}\otimes\cdots \otimes\varphi\cdot X_{q},
\]
with $\omega_{1},\ldots,\omega_{p}\in\Omega^{1}(M)$, $X_{1},
\ldots,X_{q}\in \mathfrak{X}(M)$.


Let $\pi\colon FM\rightarrow M$ be the bundle of linear frames of
$M$ and let $G$ be a closed subgroup in $GL(m;\mathbb{R})$,
$m=\dim M$. We recall (\emph{e.g.\/}, see \cite[I, pp. 57-58]{KN})
that the $G$-structures over $M$ are in bijection with the
sections of the quotient bundle $\bar{\pi}\colon F(M)/G\rightarrow
M$. In fact, the bijection $s\leftrightarrow \mathcal{B}_{s}$
between sections and $G$-structures is given by the formula
$\mathcal{B}_{s}=\left\{  u\in F(M):u\cdot G =s\left(  \pi\left(
u\right)  \right) \right\} $, where $u\cdot G$ denotes the coset
of $u\in FM$ in $FM/G$, and the inverse bijection is $s_{P}\colon
M\to FM/G$, $s_{\mathcal{B}}(x)=u(\operatorname{mod}G)$, where $u$
is any point in the fibre $\mathcal{B}_{x}$. For every
diffeomorphism $\varphi\colon M\to M^{\prime}$, let
$\tilde{\varphi}\colon FM\rightarrow FM^{\prime}$ be the $GL\left(
m;\mathbb{R}\right)  $-principal bundle isomorphism given by
(\emph{cf.\/} \cite[Chapter VI, p.\ 226]{KN}):
\[
\tilde{\varphi}\left(  X_{1},\ldots,X_{m}\right) =\left(
\varphi_{\ast }X_{1},\ldots,\varphi_{\ast}X_{m}\right) ,
\quad\left(  X_{1},\ldots ,X_{n}\right)  \in FM.
\]

Two $G$-structures $\pi\colon\mathcal{B}\rightarrow M$, $\pi\colon
\mathcal{B}^{\prime}\rightarrow M^{\prime}$ are said to be
equivalent if there exists a diffeomorphism $\varphi\colon
M\rightarrow M^{\prime}$ such that
$\tilde{\varphi}(P)=P^{\prime}$.

A diffeomorphism $\varphi\colon M\rightarrow M^{\prime}$
transforms each section $s\colon M\rightarrow F(M)/G$ of
$\bar{\pi}\colon F(M)/G\rightarrow M $ into a section of
$\bar{\pi}^{\prime}\colon F(M)/G\rightarrow M$ by setting
\[
\varphi\cdot s=\tilde{\varphi}\circ s\circ\varphi^{-1}.
\]

The $G$-structure corresponding to $\varphi\cdot s$ is said to be
obtained by transporting $\mathcal{B}_{s}$ to $M^{\prime}$ via
$\varphi$.

Two $G$-structures $\mathcal{B},\mathcal{B}^{\prime}$ are
equivalent if and only if their corresponding sections are related
by a diffeomorphism; more precisely,
$\tilde{\varphi}(\mathcal{B})=\mathcal{B}^{\prime}$ if and only if
$\varphi\cdot s_{P}=s_{P^{\prime}}$.

By passing $\tilde{\varphi}\colon FM\rightarrow FM^{\prime}$ to
the quotient modulo $G$, $\tilde{\varphi}$ induces a
diffeomorphism $\bar{\varphi}\colon FM/G\rightarrow
FM^{\prime}/G$, such that $\bar{\pi}^{\prime}\circ\bar
{\varphi}=\varphi\circ\bar{\pi}$.

We denote by $C(M)\rightarrow M$ the bundle of linear connections
on $M$. This is an affine bundle modelled over the vector bundle
$T^{\ast}(M)\otimes T^{\ast}(M)\otimes T(M)$ whose global sections
are identified to the linear connections on $M$ (\emph{cf.\/}
\cite{Kolar}).

Let $G\subseteq GL(m;\mathbb{R})$ be a Lie subgroup and let
$\mathfrak{g}$ be its Lie algebra. We denote by
$\mathfrak{g}^{(1)}$ its first prolongation; that is,
\[
\mathfrak{g}^{(1)}=\{T\in\mathrm{Hom}(\mathbb{R}^{m},\mathfrak{g})\colon
T(u)v-T(v)u=0,\forall u,v\in\mathbb{R}^{m}\},
\]
where we identify $\mathfrak{g}$ to its natural image in
$\mathfrak{gl}(m;\mathbb{R})=\mathrm{Hom}(\mathbb{R}^{m},\mathbb{R}^{m})$.

\section{The canonical connection}

Let $(M,F,P)$ be an almost biparacomplex manifold. Then the
following can be proved:

\begin{proposition}
\label{baselocal}  If $\{X_{1},\ldots,$ $X_{n}\}$ is a basis of
the distribution $T_{F}^{+}(M)$ at a point $x \in  M$, then
$\{PX_{1} , \ldots , PX_{n}\}$, $\{X_{1}+PX_{1},\ldots
,X_{n}+PX_{n}\}$ and $\{X_{1}-PX_{1},\ldots,X_{n}-PX_{n}\}$ are
basis of the distributions $T_{F}^{-}(M)$, $T_{P}^{+}(M)$, and
$T_{P}^{-}(M)$ at $x\in M$, respectively.
\end{proposition}

Now, taking Proposition \ref{baselocal} into account one can prove
that the $G$-structure determined by an almost biparacomplex
structure is given by $G=\Delta GL(n;\mathbb{R})$, where
\[
\Delta GL(n;\mathbb{R})=\left\{ \left(
\begin{array}
[c]{cc}
A & 0\\
0 & A
\end{array}
\right) \colon A \in GL(n; \mathbb{R}) \right\},
\]
is a closed subgroup of the Lie group $GL(n;\mathbb{R})$.

Then,  we can state the following:

\begin{proposition} The set of
all almost biparacomplex structures on a manifold $M$ can be
identified with the sections of the quotient bundle $FM/\Delta
GL(n;\mathbb{R})\to M$.
\end{proposition}

Now we shall define a canonical connection on an almost
biparacomplex manifold $(M,F,P)$. This connection $\nabla$ has the
following properties: (1) $\nabla$ is torsion-free iff the
distributions associated to the almost biparacomplex structure are
involutive; (2) $\nabla$ is locally flat (\emph{i.e.\/}, its
torsion and curvature tensor fields vanish) iff the $\Delta
GL(n;\mathbb{R})$-structure defined by $(F,P)$ is integrable.

\begin{lemma}
\label{equiv} Let $(F,P)$ be an almost biparacomplex structure on
a manifold $M$, let $(V_{1},V_{2},V_{3})$ be its associated
$\alpha$-structure, and let $\nabla$ be a linear connection on
$M$. The following conditions are equivalent:

\begin{enumerate}
\item [\emph{(1)}]$\nabla_{X}Y_{i}\in V_{i}$, $\forall
X\in\mathfrak{X}(M)$, $\forall Y_{i}\in V_{i}$, $i=1,2,3$.

\item[\emph{(2)}] $\nabla F=0=\nabla P$.
\end{enumerate}
\end{lemma}

One can prove the previous lemma by direct calculations taking
into account that $F$ is an isomorphism between the distributions
$V_{3}$ and $V_{4}$ and the definition of the tensor fields
$\nabla P$ and $\nabla F$ of type $(1,1)$.

\begin{theorem}
\label{teor:cruceanu}Let $(F,P)$ be an almost biparacomplex
structure on $M$ and let $(V_{1},V_{2},V_{3})$ be its associated
$\alpha$-structure. Then there exists a unique linear connection
$\nabla$ on $M$ verifying the following conditions:
\begin{enumerate}
\item[\emph{(i)}] $\nabla F=0=\nabla P$.

\item[\emph{(ii)}] $T ( X , Y ) = 0$, for all $X \in V_{1} , \, Y
\in V_{2}$, where $T$ denotes the torsion tensor field of
$\nabla$.
\end{enumerate}
Then $\nabla$ is called the \emph{canonical connection} of
$(F,P)$.
\end{theorem}

\begin{proof}
Let us assume that $\nabla$ is a linear connection satisfying the
above conditions. If $X\in V_{1},\,Y\in V_{2}$ one has
\[
0=T(X,Y)=\nabla_{X}Y-\nabla_{Y}X-[X,Y]
=\nabla_{X}Y-\nabla_{Y}X-F^{+}[X,Y]-F^{-}[X,Y],
\]
where $F^{+}$ (resp.\ $F^{-}$) denotes the projection over
$T_{F}^{+}(M)$ (resp.\ $T_{F}^{-}(M)$). On the other hand, as we
assume $\nabla$ parallelize both $F$ and $P$, then by Lemma
\ref{equiv} we have $\nabla_{X}Y\in V_{2},\nabla_{Y}X\in V_{1}$,
$X$ belonging to $V_{1}$ and $Y$ to $V_{2}$, and then one has
$\nabla_{X}Y-F^{-}[X,Y]=0,\quad-\nabla_{Y}X-F^{+}[X,Y]=0$,
\emph{i.e.\/},
\begin{equation}
\nabla_{X}Y=F^{-}[X,Y],\quad\nabla_{Y}X=F^{+}[Y,X].
\label{eq:torformula}
\end{equation}

Now we prove that this equation (\ref{eq:torformula}) determines
completely $\nabla$. Let $X,Y\in\mathfrak{X}(M)$. Then one can
decompose $X=F^{+}X+F^{-}X,\quad Y=F^{+}Y+F^{-}Y$, and then
\begin{equation}
\nabla_{X}Y=\nabla_{F^{+}X}F^{+}Y+\nabla_{F^{+}X}F^{-}Y
+\nabla_{F^{-}X}F^{+}Y+\nabla_{F^{-}X}F^{-}Y. \label{eq:madre}
\end{equation}

Taking into account the above equation (\ref{eq:torformula}) we
obtain
\begin{equation}
\nabla_{F^{+}X}F^{-}Y=F^{-}[F^{+}X,F^{-}Y],
\quad\nabla_{F^{-}X}F^{+}Y =F^{+}[F^{-}X,F^{+}Y],\label{eq:tor1}
\end{equation}
which allows us to deduce the following equations, taking into
account that $\nabla P=0$:
\[
\nabla_{F^{+}X}F^{+}Y=\nabla_{F^{+}X}P^{2}F^{+}Y
=P(\nabla_{F^{+}X}PF^{+}Y)=PF^{-}[F^{+}X,PF^{+}Y],
\]
\[
\nabla_{F^{-}X}F^{-}Y=\nabla_{F^{-}X}P^{2}F^{-}Y
=P(\nabla_{F^{-}X}PF^{-}Y)=PF^{+}[F^{-}X,PF^{-}Y],
\]

Then, by using the relations $P\circ F^{-}=F^{+}\circ P, P\circ
F^{+}=F^{-}\circ P$, we can conclude:
\begin{equation}
\nabla_{F^{+}X}F^{+}Y=F^{+}P[F^{+}X,PF^{+}Y],
\quad\nabla_{F^{-}X}F^{-}Y =F^{-}P[F^{-}X,PF^{-}Y],\label{eq:tor2}
\end{equation}
for all $X,Y\in\mathfrak{X}(M)$. Finally, joining the above
equations (\ref{eq:madre}), (\ref{eq:tor1}) and (\ref{eq:tor2}) we
obtain the general expression of the derivation law of $\nabla$:
\begin{align}
\nabla_{X}Y  & = F^{+}\Bigr([F^{-}X,F^{+}Y]+P[F^{+}X,PF^{+}Y]\Bigr
)\label{eq:conexion}\\
&+ F^{-}\Bigr([F^{+}X,F^{-}Y]+P[F^{-}X,PF^{-}Y]\Bigr),\nonumber
\end{align}
for all vector fields $X,Y$ on $M$.

Finally, we prove that the connection $\nabla$ defined by the
above equation (\ref{eq:conexion}) verifies both conditions of the
present theorem. It is an easy exercise to check that $\nabla$ is
a linear connection and that $\nabla$ verifies condition (ii),
because this condition was the starting point of our construction.

In order to prove condition i), one can prove, by Lemma
\ref{equiv}, that $\nabla$ preserves the distributions $V_{i}$,
$i=1,2,3$. Let $X$ a vector field on $M$, and let $Z_{i}\in V_{i},
i=1,2,3$. Then we have
\[
\nabla_{X}Z_{1}=F^{+}\Bigr([F^{-}X,Z_{1}]+P[F^{+}X,PZ_{1}]\Bigr)\Rightarrow
\nabla_{X}Z_{1}\in V_{1};
\]
\[
\nabla_{X}Z_{2}=F^{-}\Bigr([F^{+}X,Z_{2}]+P[F^{-}X,PZ_{2}]\Bigr)\Rightarrow
\nabla_{X}Z_{2}\in V_{2};
\]
and, finally,
\begin{align*}
\nabla_{X}Z_{3}  & =F^{+}\Bigr([F^{-}X,F^{+}Z_{3}]+P[F^{+}X,PF^{+}Z_{3}]\Bigr)\\
&+ F^{-}\Bigr([F^{+}X,F^{-}Z_{3}]+P[F^{-}X,PF^{-}Z_{3}]\Bigr).
\end{align*}
Taking into account  the equations $P\circ F^{+}=F^{-}\circ P$ and
$P\circ F^{-}=F^{+}\circ P$, we obtain
\begin{align*}
P(\nabla_{X}Z_{3})  & =F^{-}\Bigr(P[F^{-}X,F^{+}Z_{3}]+[F^{+}X,F^{-}Z_{3}]\Bigr)\\
&+ F^{+}\Bigr(P[F^{+}X,F^{-}Z_{3}]+[F^{-}X,F^{+}Z_{3}]\Bigr).
\end{align*}

Then,
$P(\nabla_{X}Z_{3})=\nabla_{X}Z_{3}\Rightarrow\nabla_{X}Z_{3}\in
V_{3}$, as wanted.
\end{proof}

\vspace{3mm}

Now we shall show that $\nabla$ measures the involutiveness of the
distributions associated to an almost biparacomplex structure and
the integrability of the $\Delta GL(n; \mathbb{R})$-structure. In
the first case, the Fr\"{o}licher-Nijenhuis tensor field of the
pair $(F,P)$ is also useful. Taking the equation $F\circ P+P\circ
F=0$ into account we have that the Nijenhuis bracket of $F$ and
$P$ is
\begin{align*}
[ F,P] (X, Y)  &  = [FX , PY] + [PX ,FY] + PF [ X , Y ] + FP [X,Y]\\
& - F[PX,Y] - F[X,PY] -P[FX,Y] -P[X ,FY]\\
&  =[FX , PY] + [PX ,FY] -F[PX,Y]\\
&  - F[X,PY] -P[FX,Y] -P[X ,FY].
\end{align*}

\begin{lemma}
Let $M$ be a manifold endowed with an almost biparacomplex structure $(F,P)$ and let $T$ be the torsion tensor
field of the  canonical connection. Then the following relations hold
\begin{eqnarray*}
[ F,P] (X,Y) = 2 PT (X , Y ), \; \forall\, X , Y \in T^{+}_{F}
(M),
\end{eqnarray*}
\begin{eqnarray*}
[ F,P] (X,Y) = -2 PT (X , Y ), \; \forall\, X , Y \in T^{-}_{F}
(M),
\end{eqnarray*}
\begin{eqnarray*}
[F,P] (X , Y) = 2 F^{-} [ X , PY] - 2 F^{+} [ PX , Y], \;
\forall\, X \in T^{+}_{F} (M), Y \in T^{-}_{F} (M).
\end{eqnarray*}
\end{lemma}

This technical result follows from an easy---but rather
long---calculation. Then, we can state

\begin{theorem}
\label{teor:fn-torsion}Let $(F,P)$ be an almost biparacomplex
structure on $M $ and let $\nabla$ be its canonical connection.
Then the following three conditions are equivalent:

\begin{enumerate}
\item [\emph{(i)}]$(F,P)$ is a biparacomplex structure.

\item[\emph{(ii)}] $[F,P]=0$.

\item[\emph{(iii)}] $T=0$, \emph{i.e.\/}, $\nabla$ is torsion-free.
\end{enumerate}
\end{theorem}

The proof of this result follows directly from the properties of
the canonical connection of $(F,P)$ and the equations of above lemma.

\begin{theorem}
\label{teor:cplana}Let $M$ be a $2n$-dimensional manifold endowed
with an almost biparacomplex structure $(F,P)$ and let $\nabla$ be
its canonical connection. Then, the following three conditions are
equivalent:

\begin{enumerate}
\item [\emph{(i)}]$T=0$ and $R=0$, $T$ and $R$ being the torsion
and curvature tensor fields of $\nabla$; \emph{i.e.\/}, the
canonical connection is locally flat.

\item[\emph{(ii)}] For every point $x\in M$ there exists an open
neighbourhood $U$ of $x$ and local coordinates
$(x_{1},\ldots,x_{n},y_{1},\ldots,y_{n})$ on $U$ such that
\begin{align*}
F\left(  \partial/\partial x_{i}\right)   & =\partial/\partial
y_{i},\quad
F\left(  \partial/\partial y_{i}\right)  =\partial/\partial x_{i},\\
\noalign{\smallskip} P\left(  \partial/\partial x_{i}\right)   &
=\partial/\partial x_{i},\quad P\left(  \partial/\partial
y_{i}\right) =-\partial/\partial y_{i},
\end{align*}
for all $i=1,\ldots,n.$

\item[\emph{(iii)}] The $\Delta GL(n;\mathbb{R})$-structure on $M$
defined by $(F,P)$ is integrable.
\end{enumerate}
\end{theorem}

\begin{proof}
The equivalence (ii)$\Leftrightarrow($iii) consists on the
adaptation of the general result about the integrability of a
$G$-structure to our case of $G=\Delta GL(n;\mathbb{R})$. Thus, we
only prove (i)$\Leftrightarrow($ii).

(i)$\Rightarrow$(ii). As $\nabla$ is locally flat, for every $x\in
M$ there exists a coordinate neighbourhood
$(U;x_{1},\ldots,x_{2n})$ such that $\nabla_{\partial/\partial
x_{i}}(\partial/\partial x_{j})=0,\quad1\leq i,j\leq2n$. Let
\[
F=\sum_{i,j=1}^{2n}f_{ji}\frac{\partial}{\partial x_{i}}\otimes
dx_{j}, P=\sum_{i,j=1}^{2n}g_{ji}\frac{\partial}{\partial
x_{i}}\otimes dx_{j}, \quad f_{ji},g_{ji}\in C^{\infty}(U)
\]
be the local expressions of $F$ and $P$ on $U$. First, we shall
show that the functions $f_{ji},g_{ji}$ are constant functions on
$U$. Let $X=\sum _{k=1}^{2n}X^{k}(\partial/\partial x_{k})$,
$X^{k}\in C^{\infty}(U)$, be any vector field on $U$. We have
\begin{align}
\nabla_{X}F  & =\sum_{i,j=1}^{2n}X(f_{ji})\frac{\partial}{\partial
x_{i}} \otimes dx_{j}+\sum_{i,j=1}^{2n}f_{ji}\nabla_{X}
\frac{\partial}{\partial x_{i}}\otimes dx_{j}
\label{eq:nablaf}\\
& +\sum_{i,j=1}^{2n}f_{ji}\frac{\partial}{\partial
x_{i}}\otimes\nabla _{X}dx_{j}.\nonumber
\end{align}
As the Christoffel symbols of $\nabla$ vanish on $U$ one has
$\nabla _{X}(\partial/\partial x_{i})=0$ and $\nabla_{X}dx_{i}=0,$
for all $i\in\{1,\ldots2n\}$, and then the above equation
(\ref{eq:nablaf}) reduces to
$\nabla_{X}F=\sum_{i,j=1}^{2n}X(f_{ji})(\partial/\partial
x_{i})\otimes dx_{j}$. Now, taking into account that $\nabla$
parallelizes $F$, we have $\nabla_{X}F=0$, thus proving
$X(f_{ji})=0$, and then $f_{ji}$ are constant functions. A similar
proof runs for $g_{ji}$.

Now we shall define a new local chart verifying condition (ii) of
the theorem. Let $\{X_{1}(x),\ldots,X_{n}(x)\}$ be a basis of
$T_{F,x}^{+}(M)$ with
$X_{i}(x)=\sum_{j=1}^{2n}\alpha_{j}^{i}\left(  \partial/\partial
x_{j}\right) _{x}$, $1\leq i\leq n$. As the local coefficients of
$F$ on the above chart $(U;x_{1},\ldots,x_{2n})$ are constant,
then the following $n$ vector fields define a local basis of
$T_{F}^{+}(M)$ on $U$:
\[
X_{i}(y)=\sum_{j=1}^{2n}\alpha_{j}^{i}\left(  \partial/\partial
x_{j}\right) _{y},\quad\forall y\in U,\quad1\leq i\leq n,
\]

Let us consider the vector fields on $U$ defined by
$U_{i}=X_{i}+PX_{i},\quad V_{i}=X_{i}-PX_{i},\quad1\leq i\leq n$.
By Proposition \ref{baselocal},
$\{U_{1},\ldots,U_{n},V_{1},\ldots,V_{n}\}$ is a local basis of
$TM$ at $U$. One easily check that
$[U_{i},U_{j}]=0,\,[U_{i},V_{j}]=0,\,[V_{i},V_{j}]=0,\,1\leq
i,j\leq n$, taking into account that $F$ and $P$ have constant
coefficients on $(U;x_{1},\ldots,x_{2n})$. Then there exist
coordinates $(u_{1},\ldots,u_{n},v_{1},\ldots,v_{n})$ on $U$ such
that $U_{i}=\partial/\partial u_{i}$, $V_{i}=\partial/\partial
v_{i}$, $1\leq i\leq n$. Then we have
\[
F\left(  \partial/\partial u_{i}\right)
=F(X_{i}+PX_{i})=X_{i}+FPX_{i} =X_{i}-PX_{i}=\partial/\partial
v_{i},
\]
\[
F\left(  \partial/\partial v_{i}\right)
=F(X_{i}-PX_{i})=X_{i}-FPX_{i} =X_{i}+PX_{i}=\partial/\partial
u_{i},
\]
\[
P\left(  \partial/\partial u_{i}\right) =P(X_{i}+PX_{i})=PX_{i}
+X_{i}=\partial/\partial u_{i},
\]
\[
P\left(  \partial/\partial v_{i}\right) =P(X_{i}-PX_{i})=PX_{i}
-X_{i}=-\partial/\partial v_{i},
\]
for all $1\leq i\leq n$, as wanted.

(ii)$\Rightarrow$(i). The existence of such a local chart implies that the almost paracomplex structures defined
by $F$ and $P$ are in fact paracomplex, and then the associated distributions are involutive (see \cite[Prop. 1.2]{KK}). Then the structure is biparacomplex and, by the above Theorem \ref{teor:fn-torsion}, the canonical
connection is torsion-free, thus proving
\begin{equation}
T\left(  \partial/\partial x_{i},\partial/\partial y_{j}\right)
=\nabla_{\partial/\partial x_{i}}(\partial/\partial y_{j})-\nabla
_{\partial/\partial y_{j}}(\partial/\partial x_{i})=0.
\label{eq:torsioncero}
\end{equation}
We only need to prove that $R$ also vanishes. As the structure is
biparacomplex, we know that $\nabla F=0=\nabla P$. By Lemma
\ref{equiv} we can deduce that $\nabla_{\partial/\partial x_{i}}
(\partial/\partial y_{j})\in T_{P}^{-}(M)$ and that
$\nabla_{\partial/\partial y_{j}}(\partial/\partial x_{i})\in
T_{P}^{+}(M)$. Taking into account that $TM=T_{P}^{+}(M)\oplus
T_{P}^{-}(M)$, from the above equation (\ref{eq:torsioncero}) we
obtain
\[
\nabla_{\partial/\partial x_{i}}(\partial/\partial
y_{j})=0,\,\nabla _{\partial/\partial y_{j}}(\partial/\partial
x_{i})=0,\quad1\leq i,j\leq n.
\]

Moreover, as $F^{2}=Id$ and $\nabla F=0$, for all
$i,j\in\{1,\ldots,n\}$ we can deduce
\[
\nabla_{\partial/\partial x_{i}}(\partial/\partial x_{j})=F\left(
\nabla_{\partial/\partial x_{i}}F\left(
\partial/\partial x_{j}\right) \right)  =F\left(  \nabla_{\partial/\partial x_{i}}\left(  \partial/\partial
y_{j}\right)  \right)  =0,
\]
\[
\nabla_{\partial/\partial y_{i}}(\partial/\partial y_{j})=F\left(
\nabla_{\partial/\partial y_{i}}F\left(
\partial/\partial y_{j}\right) \right)  =F\left(  \nabla_{\partial/\partial y_{i}}\left(  \partial/\partial
x_{j}\right)  \right)  =0,
\]
for all $i,j=1,\ldots,n$. So we have proved that all the
Christoffel symbols vanish, and hence $\nabla$ is locally flat.
\end{proof}

\medskip

We end this study obtaining the expression of the canonical
connection on an adapted local frame. This result will be useful
in order to compare the canonical connection with other
connections. Let $M$ be a manifold endowed with an almost
biparacomplex structure $(F,P)$. Let $\nabla$ be its canonical
connection. Let $U$ be an open subset of $M$ and let
$\{X_{1},\ldots,X_{n},Y_{1},\ldots,Y_{n}\}$ be a local frame on
$U$ adapted to $(F,P)$. Let us denote by
$\{\omega_{1},\ldots,\omega_{n},\eta_{1},\ldots,\eta_{n}\}$ its
dual coframe. Then, $\nabla$ is determined on $U$ by
$\nabla_{X_{h}}X_{a}$, $\nabla_{Y_{h}}X_{a}$, $a,h=1,\ldots,n$,
because $\nabla_{X_{h}}Y_{a}=P\nabla_{X_{h}}X_{a}$,
$\nabla_{Y_{h}}Y_{a}=P\nabla_{Y_{h}}X_{a}$, $a,h=1,\ldots,n$.

As $\nabla_{X_{h}}X_{a},\nabla_{Y_{h}}X_{a}\in T_{F}^{+}(M)$ for
all $h,a=1,\ldots,n$, we can write
\[
\nabla_{X_{h}}X_{a}=\sum_{i=1}^{n}\Gamma_{ha}^{i}X_{i},
\quad\nabla_{Y_{h}}X_{a}=\sum_{i=1}^{n}{\bar{\Gamma}}_{ha}^{i}X_{i},
\quad a,h=1,\ldots,n,
\]
and from $T(X_{h},Y_{a})=0$ we obtain
$\sum_{i=1}^{n}(\Gamma_{ha}^{i}Y_{i}-{\bar{\Gamma}}_{ah}^{i}X_{i})
=[X_{h},Y_{a}]$. Hence
\[
\sum_{i=1}^{n}\Gamma_{ha}^{i}Y_{i}
=\sum_{i=1}^{n}\eta_{i}([X_{h},Y_{a}])Y_{i},
\quad-\sum_{i=1}^{n}{\bar{\Gamma}}_{ah}^{i}X_{i}
=\sum_{i=1}^{n}\omega_{i}([X_{h},Y_{a}])X_{i},
\]
and finally,

\begin{equation}
\Gamma_{ha}^{i}=\eta_{i}([X_{h},Y_{a}]),
\quad{\bar{\Gamma}}_{ah}^{i} =\omega_{i}([Y_{a},X_{h}]),\quad
a,h,i=1,\ldots,n. \label{eq:canonicalocal}
\end{equation}

\subsection{Functorial connections}

Now we shall prove that the canonical connection of an almost
biparacomplex structure is a functorial connection. Roughly
speaking, a functorial connection associated to a $G$-structure is
a family of reducible connections, one for each concrete
$G$-structure, which is natural with respect to the isomorphisms
of the $G$-structure. We shall also show that in our case
$G=\Delta GL(n;\mathbb{R})$ there exists, at least, another
functorial connection, which will be called the \emph{well-adapted
connection}. Both connections, the canonical and the well-adapted,
coincide iff the manifold is a biparacomplex manifold.

\vspace{3mm} Now, we shall present general results concerning
functorial connections.
\begin{definition}

\label{teor:deffuntorial}
 Let $G$ be a Lie subgroup of
$GL(m;\mathbb{R})$, with $m=\dim M$. A \emph{functorial
connection} attached to the $G$-structures over $M$ is a presheaf
morphism$\;s\mapsto\nabla\left(  s\right)  $, from the sheaf of
sections of $FM/G\rightarrow M$ into the presheaf of sections of
$C(M)\rightarrow M$, such that:

\begin{enumerate}
\item [(1)] If $s\colon U\rightarrow FM$ is a section on an open subset $U\subseteq M$, then $\nabla \left(
s\right) $ is a linear connection on $U$ adapted to the $G$-structure defined by $s$; \emph{i.e.\/},
$\nabla\left( s\right) $ is reducible to the subbundle $\mathcal{B}_{s}$ (see \S\ref{notandprel}).

\item[(2)] If $\varphi\colon U\rightarrow U^{\prime}$ is a
diffeomorphism, then $\nabla\left(  \varphi\cdot s\right)
=\varphi\cdot\nabla\left( s\right)  $, where
$\varphi\cdot\nabla\left(  s\right)  $ denotes the direct image of
$\nabla\left(  s\right)  $ by\emph{\ }$\varphi$ (\emph{cf.\/}
\cite[II. Prop. 6.1]{KN}).

\item[(3)] There exists a non-negative integer $r$ such that
$\nabla$ factors smoothly through $J^{r}\left( FM/G\right)  $.
\end{enumerate}
\end{definition}

The third item above means that the value of the section
$\nabla\left( s\right)  $ of $C(U)\rightarrow U$ at a point $x\in
U$ depends only on $j_{x}^{r}s$, and that the induced map
\begin{align*}
\nabla^{r}\colon J^{r}\left(  FM/G\right)   & \rightarrow C(M)\\
\nabla^{r}\left(  j_{x}^{r}s\right)   & =\nabla\left(  s\right)
(x)
\end{align*}
is differentiable. Also note that $\nabla=\nabla^{r}\circ j^{r}$.
This item can be substituted by apparently less restrictive
conditions; for example, by only imposing that $\nabla\left(
s\right)  (x)$ depends on the germ of $s$ at $x$, not on the
$r$-jet of the section. Nevertheless, standard techniques working
in a very general setting (see \cite[Chapter V]{Kolar}) readily
shows the equivalence of both conditions.

\begin{theorem}
\label{teor:metodo} \emph{\cite[Th. 1.1]{Valdes} } The following
conditions are equivalent:
\begin{enumerate}
\item[\emph{(i)}] For every $G$-structure $\mathcal{B}\rightarrow
M$ there exists a unique connection $\nabla^{\prime}$ adapted to
the $G$-structure such that, for every endomorphism $S$ given by a
section of the adjoint bundle $\mathrm{ad}\mathcal{B}$ and every
vector field $X\in\mathfrak{X}(M)$, one has $\mathrm{trace}(S\circ
i_{X}\circ T^{\prime})=0$, where $T^{\prime}$ denotes the torsion
of $\nabla^{\prime}$. Moreover, this connection only depends on
the first contact of the $G$-structure.

\item[\emph{(ii)}] If
$L\in\mathrm{Hom}(\mathbb{R}^{n},\mathfrak{g})$ verifies
$i_{v}\circ\mathrm{Alt}(L)\in\mathfrak g^{\perp}$ for every
$v\in^{n}$ then $L=0$, where $\mathfrak g$ is the Lie algebra of
the group $G$, $\mathfrak g^{\perp}$ is the orthogonal subspace of
$\mathfrak g$ in $\mathfrak{gl}(n;\mathbb{R})$ respect to the
Killing-Cartan metric, and $\mathrm{Alt}(L)(u,v)=L(u)v-L(v)u$,
$\forall u,v\in\mathbb{R}^{n}$.
\end{enumerate}
If there exists, we shall call this connection the
\emph{well-adapted connection} to the $G$-structure
$\pi\colon\mathcal{B}\rightarrow M$.
\end{theorem}

The well-adapted connection is a functorial connection and
measures the integrability of the $G$-structure in the sense that
the $G$-structure is integrable iff the well-adapted connection is
locally flat (\emph{cf.\/} \cite[Th. 2.3]{Valdes}). The above
Theorem \ref{teor:metodo} gives us an explicit way of obtaining a
functorial connection. Moreover, one has

\begin{theorem}
\emph{\cite[Th. 2.1]{Valdes} } Let $G$ be a Lie group with Lie
algebra $\mathfrak{g}$. If $\mathfrak{g}^{(1)}=0$ and
$\mathfrak{g}$ is invariant under transposition, then condition
ii) of \emph{Theorem \ref{teor:metodo}} holds.
\end{theorem}

Then, we have finished the study of the general framework. Now, we
specialize to the case $G=\Delta GL(n;\mathbb{R})$, showing that
(1) the canonical connection is a functorial connection; (2) there
exists the well-adapted connection; and (3) the canonical and the
well-adapted connections coincide iff the manifold is
biparacomplex.

\vspace{3mm}

The following result will be useful to prove that the canonical
connection is functorial. In fact, this is a property stronger
than that of a functorial connection.

\begin{theorem}
\label{teor:functorial}
 Let $M,M^{\prime}$ be two manifolds endowed with almost biparacomplex structures
$(F,P),(F^{\prime},P^{\prime})$, respectively and let
$\varphi\colon M\rightarrow M^{\prime}$ be a diffeomorphism
between both structures, \emph{i.e.\/},
\begin{equation}
\varphi_{\ast}\circ F
=F^{\prime}\circ\varphi_{\ast},\quad\varphi_{\ast}\circ P
=P^{\prime}\circ\varphi_{\ast}. \label{eq:fpvarphi}
\end{equation}
Let $\nabla,\nabla^{\prime}$ be the canonical connections of
$(F,P),(F^{\prime },P^{\prime})$, respectively. Then
$\nabla^{\prime}$ is the direct image of $\nabla$ via $\varphi$.
\end{theorem}
\begin{proof}
It is a direct consequence of equation (\ref{eq:conexion}), which states $\nabla$ in terms of $F$ and $P$.
\end{proof}

\begin{theorem}
\label{canonicaesfuntorial} The assignment of the canonical
connection of $(F_{\sigma},P_{\sigma})$ to each section $\sigma$
of the fiber bundle $F(M)/\Delta GL(n;\mathbb{R})\rightarrow M$ is
a functorial connection, $(F_{\sigma},P_{\sigma})$ being the
almost biparacomplex structure associated to the section $\sigma$.
\end{theorem}

\begin{proof}
Obviously $\nabla$ is an adapted connection and then it satisfies
the first condition in Definition \ref{teor:deffuntorial}. The
above Theorem \ref{teor:functorial} shows that the second
condition is also satisfied. Finally, from equation
(\ref{eq:canonicalocal}) we obtain that $\nabla(x)$ depends only
on $f_{ij}(x)$, $g_{ij}(x)$, $\partial f_{ij}/\partial x_{k}(x)$,
$\partial g_{ij}/\partial x_{k}(x)$, $i,k=1,\ldots,2n$, $1\leq
j\leq n$, and then, taking $r=1$, the third condition is
satisfied.
\end{proof}

\vspace{3mm}

The following result proves the existence of the well-adapted
connection.

\begin{proposition}
\label{existewelladapted} The first prolongation of the Lie
algebra ${\Delta_{\ast}}\mathfrak{gl}(n;\mathbb{R})$ of the Lie
group $\Delta GL(n;\mathbb{R})$ vanishes and
${\Delta_{\ast}\mathfrak{gl}(n; \mathbb{R})}$ is invariant under
transposition. Then, there exists the well-adapted connection.
\end{proposition}

The proof of this result is an exercise of Linear Algebra and
therefore is omitted.

\medskip

Now we shall obtain the local expression of the well-adapted
connection $\nabla^{\prime}$ of an almost biparacomplex structure
$(F,P)$. We are looking for comparing with the expression of the
canonical connection obtained in equation
(\ref{eq:canonicalocal}). The first step consists on determining
the well-adapted connection by means of Theorem \ref{teor:metodo},
so we must obtain information about the adjoint bundle.

Let $\mathcal{B}\rightarrow M$ be the $\Delta
GL(n;\mathbb{R})$-structure defined by $(F,P)$. We denote by
$\mathrm{ad}\mathcal{B}$ the associated fiber bundle to
$\mathcal{B}$ via the adjoint representation of $\Delta
GL(n;\mathbb{R})$ on $\Delta_{\ast}GL(n;\mathbb{R})$,
\emph{i.e.\/},
$\mathrm{ad}\mathcal{B}=(\mathcal{B}\times{\Delta_{\ast}\mathfrak{
gl}(n; \mathbb{R})})/ \Delta GL(n;\mathbb{R})$, where the action
of $\Delta GL(n;\mathbb{R})$ on
$\mathcal{B}\times{\Delta_{\ast}}\mathfrak{gl}{(}n;\mathbb{R})$ is
given by
\[
\left(  u,\left(
\begin{array}
[c]{cc}
A & 0\\
0 & A
\end{array}
\right)  \right)  \cdot\left(
\begin{array}
[c]{cc}
B & 0\\
0 & B
\end{array}
\right)  =\left(  u\cdot\left(
\begin{array}
[c]{cc}
B & 0\\
0 & B
\end{array}
\right)  ,\right.  \left.  \left(
\begin{array}
[c]{cc}
B^{-1}AB & 0\\
0 & B^{-1}AB
\end{array}
\right)  \right)  ,
\]
where

\[
u\in\mathcal{B},\quad\left(
\begin{array}
[c]{cc}
B & 0\\
0 & B
\end{array}
\right)  \in\Delta GL(n;\mathbb{R}),\quad\left(
\begin{array}
[c]{cc}
A & 0\\
0 & A
\end{array}
\right)  \in{\Delta_{\ast}}\mathfrak{gl}(n;\mathbb{R}).
\]

One observes that the adjoint bundle $\mathrm{ad}\mathcal{B}$ is a
subbundle of the bundle of endomorphisms of the tangent bundle,
$\mathrm{ad}\mathcal{B}\subset\mathrm{End}\,(TM)\cong\mathcal{T}_{1}^{1}(M)$,
because
${\Delta_{\ast}}\mathfrak{gl}(n;\mathbb{R})\subset\mathfrak{gl}(2n;\mathbb{R})
\equiv\mathrm{End}(\mathbb{R}^{2n})$. The following result
characterizes the endomorphisms of the tangent bundle which
belongs to the adjoint bundle. Let us introduce some notations:
let $x\in M$ and let $U$ an open neighbourhood of $x$. Let
$\{X_{1},\ldots,X_{n},Y_{1},\ldots,Y_{n}\}$ be a local frame
adapted to $(F,P) $ on $U$. Then, an element of
$(\mathrm{ad}\mathcal{B})_{x}$ is, by definition, an endomorphism
$S\colon T_{x}M\rightarrow T_{x}M$
\begin{equation}
S(X_{j})=\sum_{i=1}^{n}a_{ij}X_{i},\quad S(Y_{j})
=\sum_{i=1}^{n}a_{ij}Y_{i},\quad j=1,\ldots n,\label{endomorfismo}
\end{equation}
where $A=(a_{ij})$ is any matrix of $\mathfrak{gl}(n;\mathbb{R})$.

\begin{lemma}
An endomorphism $S$ of the tangent bundle can be written as in
\emph{(\ref{endomorfismo})}, for every point $x\in M$, iff $S$
commutes with $F$ and $P$.
\end{lemma}

\begin{proof}
Let
\[
\left(
\begin{array}
[c]{cc}
A & B\\
C & D
\end{array}
\right)  ,\;A,B,C,D\in\mathfrak{gl}(n;\mathbb{R}),
\]
be the matrix of $S$ respect to the basis
$\{X_{1},\ldots,X_{n},Y_{1},\ldots,Y_{n}\}$. Respect to such
basis, the expressions of $F$ and $P$ are:
\[
\left(
\begin{array}
[c]{rr}
I_{n} & 0\\
0 & -I_{n}
\end{array}
\right)  ,\quad\left(
\begin{array}
[c]{ll}
0 & I_{n}\\
I_{n} & 0
\end{array}
\right)  .
\]
Then $F\circ S=S\circ F$ and $P\circ S=S\circ P$ iff $B=C=0$ and
$D=A$. \emph{i.e.\/},
$S\in{\Delta_{\ast}}\mathfrak{gl}{(n;\mathbb{R})}$.
\end{proof}

\medskip

Then we can deduce:

\begin{proposition}
The sections of the adjoint bundle of a
$\Delta GL(n;\mathbb{R})$-structure over $M$ defined by an almost
biparacomplex structure $(F,P)$ on $M$ are the endomorphisms of
$TM$ which commute with $F$ and $P$.
\end{proposition}

Condition (i) of Theorem \ref{teor:metodo} is given by
$\mathrm{trace}(S\circ i_{X}\circ T^{\prime})=0 ,\,\forall
S\in\Gamma(\mathrm{ad}\mathcal{B})$, $\forall
X\in\mathfrak{X}(M)$, $T^{\prime}$ being the torsion tensor of the
well-adapted connection $\nabla^{\prime}$. Taking into account the
definition of the trace, one has:
\begin{align*}
0  & =\sum_{i=1}^{n}\omega_{i}((S\circ i_{X}\circ
T^{\prime})(X_{i}))
+\sum_{i=1}^{n}\eta_{i}((S\circ i_{X}\circ T^{\prime})(Y_{i}))\\
\noalign{\smallskip}  &
=\sum_{i=1}^{n}\omega_{i}(S(T^{\prime})(X,X_{i})))
+\sum_{i=1}^{n}\eta_{i}(S(T^{\prime}(X,Y_{i}))),
\end{align*}
where $T^{\prime}$ denotes the torsion tensor of $\nabla^{\prime}$
and $\{\omega_{1},\ldots,\omega_{n},\eta_{1},\ldots,\eta_{n}\}$
denotes the dual coframe of
$\{X_{1},\ldots,X_{n},Y_{1},\ldots,Y_{n}\}$ on $U$. We can locally
determine $\nabla^{\prime}$ choosing a local basis of sections of
the adjoint bundle $\mathrm{ad}\mathcal{B}$ and a local basis of
$\mathfrak{X}(M) $. Let us take the family $\{\omega_{b}\otimes
X_{a}+\eta_{b}\otimes Y_{a}:\,a,b=1,\ldots , n\}$ as local basis
of $\Gamma(\mathrm{ad}\mathcal{B})$ and let us take the local
adapted frame as a local basis of $\mathfrak{X}(U)$. Then, for
$S=\omega_{b}\otimes X_{a}+\eta_{b}\otimes Y_{a}$ and $X=X_{h}$,
we have
\[
\sum_{i=1}^{n}\omega_{i}((\omega_{b}\otimes X_{a}+\eta_{b}\otimes
Y_{a})(T^{\prime}(X_{h},X_{i}))+\sum_{i=1}^{n}\eta_{i}((\omega_{b}\otimes
X_{a}+\eta_{b}\otimes Y_{a})(T^{\prime}(X_{h},Y_{i}))=0,
\]
and, consequently,
\begin{equation}
\omega_{b}(T^{\prime}(X_{h},X_{a}))+\eta_{b}(T^{\prime}(X_{h},
Y_{a}))=0.\label{xh}
\end{equation}
On the other hand, for $X=Y_{h}$ and $S=\omega_{b}\otimes
X_{a}+\eta _{b}\otimes Y_{a}$, we have
\[
\sum_{i=1}^{n}\omega_{i}((\omega_{b}\otimes X_{a}+\eta_{b}\otimes
Y_{a})(T^{\prime}(Y_{h},X_{i}))+\sum_{i=1}^{n}\eta_{i}((\omega_{b}\otimes
X_{a}+\eta_{b}\otimes Y_{a})(T^{\prime}(Y_{h},Y_{i}))=0,
\]
and then,
\begin{equation}
\omega_{b}(T^{\prime}(Y_{h},X_{a}))+\eta_{b}(T^{\prime}(Y_{h},Y_{a}))=0.
\label{yh}
\end{equation}
As $\nabla^{\prime}$ is a functorial connection, $\nabla^{\prime}$
parallelizes $F$ and $P$ and preserves the distributions
associated to the eigenvalues $\pm1$ of $F$ and $P$ (see Lemma
\ref{equiv}). Then,
$\nabla_{X_{h}}^{\prime}X_{a},\nabla_{Y_{h}}^{\prime}X_{a}\in
T_{F}^{+}(M)$,
$\,\nabla_{X_{h}}^{\prime}Y_{a},\nabla_{Y_{h}}^{\prime}Y_{a}\in
T_{F}^{-}(M)$. Moreover, one observes that $\nabla^{\prime}$ is
completely determined by the values
$\nabla_{X_{h}}^{\prime}X_{a},\nabla_{Y_{h}}^{\prime}X_{a}$,
$h,a=1,\ldots , n$. Let us denote
\[
\nabla_{X_{h}}^{\prime}X_{a}
=\sum_{i=1}^{n}\Gamma_{ha}^{^{\prime}i}X_{i} \quad,
\nabla_{Y_{h}}^{\prime}X_{a}=\sum_{i=1}^{n}{\bar{\Gamma}}_{ha}^{^{\prime}i}X_{i}.
\]
Next, we show that the equations (\ref{xh}) and (\ref{yh}) allow us to obtain the Christoffel symbols
$\Gamma_{ha}^{^{\prime}i},{\bar{\Gamma}}_{ha}^{^{\prime}i}$, $h,a,i=1,\ldots,n$, on $U$. From equation
(\ref{xh}), we have
\begin{equation}
0=2\Gamma_{ha}^{^{\prime}b} -\Gamma_{ah}^{^{\prime}b}
=\omega_{b}([X_{h},X_{a}]) +\eta_{b}([X_{h},Y_{a}]),\label{xha}
\end{equation}
and permuting the indices $a$ and $h$, we obtain
\begin{equation}
0=2\Gamma_{ah}^{^{\prime}b}-\Gamma_{ha}^{^{\prime}b}
=\omega_{b}([X_{a},X_{h}])+\eta_{b}([X_{a},Y_{h}]). \label{xah}
\end{equation}
From equations (\ref{xha}) and (\ref{xah}) above we conclude
\begin{equation}
\Gamma_{ah}^{^{\prime}b}
=\frac{1}{3}\Bigr(\omega_{b}([X_{a},X_{h}]) +2\eta
_{b}([X_{a},Y_{h}])+\eta_{b}([X_{h},Y_{a}])\Bigr) . \label{gammas}
\end{equation}

Moreover, by using equation (\ref{yh}), and by a similar argument
to the one above, we can obtain the remaining Christoffel symbols:

\begin{equation}
{\bar\Gamma}_{ah}^{^{\prime}b} = \frac{1}{3} \Bigr( \omega_{b} (
[Y_{h}, X_{a} ]) + 2 \omega_{b} ( [Y_{a} , X_{h} ]) + \eta_{b} (
[Y_{a},Y_{h} ])\Bigr).\label{gammasbarra}
\end{equation}

\begin{remark}
Let $(F,P)$ be an almost biparacomplex structure on a manifold $M$ and let $\nabla$ (resp. $\nabla^{\prime}$) be
its canonical (resp. well-adapted) connection. In general, these connections do not coincide. The proof is very
easy: one only must compare equations (\ref{eq:canonicalocal}) of $\nabla$ with equations (\ref{gammas}) and
(\ref{gammasbarra}) of $\nabla^{\prime}$.
\end{remark}

Moreover, we can obtain the expression of the tensor field of type
$(1,2)$ given by the difference of both connections:
$A=\nabla-\nabla^{\prime}$.

\begin{theorem}
Let $(F,P)$ be an almost biparacomplex structure on a manifold $M$
and let $\nabla$ (resp. $\nabla^{\prime}$) be its canonical (resp.
well-adapted) connection. Let $A$ be the tensor defined by $A(X,Y)
= \nabla_{X} Y - \nabla^{\prime}_{X} Y , \, \forall X , Y
\in\mathfrak{X}(M)$. Then,
\begin{align*}
A(X,Y)  & =\frac{1}{3}\bigg(F^{+}T(F^{+}X,F^{+}Y)+PF^{+}T(F^{+}X,PF^{-}Y)\\
\noalign{\smallskip}  &
+PF^{-}T(F^{-}X,PF^{+}Y)+F^{-}T(F^{-}X,F^{-}Y)\bigg).
\end{align*}
for all vector fields $X,Y$ on $M$, where $T$ denotes the torsion
tensor of the canonical connection.
\end{theorem}

\begin{proof}
From equations (\ref{eq:canonicalocal}) and (\ref{gammas}) we
deduce
\begin{align*}
A(X_{h},X_{a})  & =\nabla_{X_{h}}X_{a}-\nabla_{X_{h}}^{\prime}X_{a}\\
\noalign{\smallskip}  & =\frac{1}{3}\bigg(PF^{-}[X_{h},PX_{a}]
-F^{+}[X_{h},X_{a}]+PF^{-}[PX_{h},X_{a}]\bigg)\\
\noalign{\smallskip}  & =\frac{1}{3}F^{+}T(X_{h},X_{a}), \quad
\forall h,a=1,\ldots,n.
\end{align*}

Then
\begin{align}
A(X^{\prime},X^{\prime\prime})  & =\frac{1}{3}
\bigg(PF^{-}[X^{\prime },PX^{\prime\prime}]
-F^{+}[X^{\prime},X^{\prime\prime}]+PF^{-}[PX^{\prime },
X^{\prime\prime}]\bigg)\nonumber\\
\noalign{\smallskip}  &
=\frac{1}{3}F^{+}T(X^{\prime},X^{\prime\prime }), \label{eq:t11}
\end{align}
for all $X^{\prime},X^{\prime\prime}\in T_{F}^{+}(M)$. Moreover,
we have
\[
A(X_{h},Y_{a})=\nabla_{X_{h}}Y_{a}-\nabla_{X_{h}}^{\prime}Y_{a}
=P(\nabla
_{X_{h}}X_{a}-\nabla_{X_{h}}^{\prime}X_{a})=PA(X_{h},PY_{a}),
\]
for all $h,a=1,\ldots,n$. Then
\begin{align}
A(X^{\prime},Y^{\prime})  & =\frac{1}{3}
\bigg(F^{-}[X^{\prime},Y^{\prime }]
-PF^{+}[X^{\prime},PY^{\prime}]+F^{-}[PX^{\prime},PY^{\prime}]\bigg
)
\nonumber\\
\noalign{\smallskip}  & =\frac{1}{3}PF^{+}T(X^{\prime},PY^{\prime
}), \label{eq:t12}
\end{align}
for all $X^{\prime}\in T_{F}^{+}(M),Y^{\prime}\in T_{F}^{-}(M)$.

Similar arguments allow us to obtain the following equalities:
\begin{align}
A(Y^{\prime},X^{\prime}) &
=\frac{1}{3}\bigg(F^{+}[Y^{\prime},X^{\prime
}]+F^{+}[PY^{\prime},PX^{\prime}]-PF^{-}[Y^{\prime},PX^{\prime}]\bigg
)\nonumber\\
\noalign{\smallskip} & =\frac{1}{3}PF^{-}T(Y^{\prime},PX^{\prime
}),\label{eq:t21}
\end{align}
for all $X^{\prime}\in T_{F}^{+}(M),Y^{\prime}\in T_{F}^{-}(M)$;
\begin{align}
A(Y^{\prime},Y^{\prime\prime}) & =\frac{1}{3}
\bigg(PF^{+}[Y^{\prime },PY^{\prime\prime}]
+PF^{+}[PY^{\prime},Y^{\prime\prime}]-F^{-}[Y^{\prime },
Y^{\prime\prime}]\bigg)\nonumber\\
\noalign{\smallskip} &
=\frac{1}{3}F^{-}T(Y^{\prime},Y^{\prime\prime }), \label{eq:t22}
\end{align}
for all $Y^{\prime},Y^{\prime\prime}\in T_{F}^{-}(M)$.

From the above equations (\ref{eq:t11}), (\ref{eq:t12}), (\ref{eq:t21}) and (\ref{eq:t22}) we obtain two
expressions for the tensor $A$; in the first one $A$ is refereed to the tensors $F,P,F^{+}$ and $F^{-}$:
\begin{align*}
A(X,Y)  & =\frac{1}{3}\bigg(PF^{-}[FX,PF^{+}Y]
+PF^{+}[PF^{+}X,F^{+}Y]-F^{+}[FX,F^{+}Y]\\
\noalign{\smallskip}  & +F^{+}[PF^{-}X,PF^{+}Y]
+F^{-}[FX,F^{-}Y]+F^{-}[PF^{+}X,PF^{-}Y]\\
\noalign{\smallskip}  &
-PF^{+}[FX,PF^{-}Y]+PF^{+}[PF^{-}X,F^{-}Y]\bigg),
\end{align*}
and in the second one, $A$ is refereed to the torsion $T$ of the
canonical connection:
\begin{align}
A(X,Y)  &
=\frac{1}{3}\bigg(F^{+}T(F^{+}X,F^{+}Y)+PF^{+}T(F^{+}X,PF^{-}Y)
\label{eq:diferenciatorsion}\\
\noalign{\smallskip}  &
+PF^{-}T(F^{-}X,PF^{+}Y)+F^{-}T(F^{-}X,F^{-}Y)\bigg ). \nonumber
\end{align}
\end{proof}

\vspace{3mm}

As a direct consequence of this result we obtain

\begin{theorem}
\label{coinciden}
 Let $(F,P)$ be an almost biparacomplex structure on a manifold $M$ and let $\nabla$ (resp.
$\nabla^{\prime}$) be its canonical (resp. well-adapted)
connection. If $(F,P)$ is a biparacomplex structure then $\nabla
=\nabla^{\prime}$.
\end{theorem}

\begin{proof}
By Theorem \ref{teor:fn-torsion}, if $(F,P)$ is biparacomplex, then $\nabla$ is torsion-free and hence, from
equation (\ref{eq:diferenciatorsion}), we obtain $A=0$, thus proving $\nabla=\nabla^{\prime}$.
\end{proof}

\section{Examples and Final Remarks}

We finish the paper showing some examples and point out some
aspects of the theory.

\subsection{Equivalence of $\Delta GL(n;\mathbb R )$-structures}

The equation (\ref{eq:fpvarphi}) of Theorem \ref{teor:functorial}
is the definition of the equivalence of the almost biparacomplex
structures $(F,P)$ and $(F',P')$ over the manifolds $M$ and $M'$:
we say that $(F,P)$ and $(F',P')$ are {\em equivalent} if there
exists a diffeomorphism $\varphi \colon M\to M'$ such that $
\varphi_{\ast}\circ F
=F^{\prime}\circ\varphi_{\ast},\quad\varphi_{\ast}\circ P
=P^{\prime}\circ\varphi_{\ast}$.

This condition is equivalent to the classical definition of
$G$-structures in the case of the Lie group $\Delta GL(n;\mathbb R
)$.

\begin{proposition} \emph{\cite{San2} } Let $\pi \colon {\cal B} \to M$ y $\pi '
\colon {\cal B}' \to M'$ be the $\Delta GL(n;\mathbb R
)$-structures over the $n$-dimensional manifolds $M$ and $M'$
defined by the almost biparacomplex structures  $(F,P)$ and
$(F',P')$ respectively. Then we have that $\cal B$ and $\cal B'$
are equivalent iff the structures $(F,P)$ and $(F',P')$
are equivalent.
\label{teor:equivalenciag}
\end{proposition}

Theorem \ref{teor:functorial} establishes that if the almost
biaparacomplex structures $(F,P)$ and $(F',P')$ are equivalent
then the canonical connections $\nabla$ and $\nabla'$ are
equivalent; {\em i.e.}, $\nabla'$ is the direct image of $\nabla$
via $\varphi$. This result is the key that allow us to
characterize the equivalence problem of $\Delta GL(n;\mathbb R
)$-structures in terms of the linear connections in the analytic
case. We have that

\begin{theorem} \emph{\cite{San2} }
Let  $(F,P), (F ' ,P ' )$ be two analytic almost biparacomplex
structures over the analytic manifolds $M , M '$ respectively, and
let $ \varphi\, \colon M \to M '$ be an analytic diffeomporphism
verifying that:
\newcounter{Item501}
\begin{list}{ \roman{Item501}$\null\,)$} {\usecounter{Item501} }

\item  there exists a point  $x_0 \in M$such that $(F,P)$ and $ (F
' ,P ' )$ are equivalent by $ \varphi$;
     \[ F_{\varphi (x_0) } ' \circ (\varphi_* )_{x_0} = (\varphi_* )_{x_0} \circ F_{x_0}  , \quad P_{\varphi (x_0) } ' \circ (\varphi_* )_{x_0} = (\varphi_* )_{x_0} \circ P_{x_0},\]
       \item the canonical connections $\nabla$ and $\nabla'$ are equivalent by $\varphi$.
\end{list}
In these conditions,  there exists an open neighbourhood $U$ of $x_0$ such that  $(F,P)$ and $  (F ' ,P ' )$ are equivalent $\varphi$ on $U$; {\em i. e.},
\[F_{\varphi (x) } ' \circ (\varphi_* )_{x} = (\varphi_* )_{x} \circ F_{x}  , \quad P_{\varphi (x) } ' \circ (\varphi_* )_{x} = (\varphi_* )_{x} \circ P_{x}, \quad \forall x \in U . \]

Moreover, if $M$ is a  connected manifold then the open $U$
coincides with $M$.
\label{teor:equivalencia}
\end{theorem}

We leave out the proofs of Proposition \ref{teor:equivalenciag} and Theorem \ref{teor:equivalencia} of this work, which can be found in \cite{San2}.
\vspace{3mm}

The main objective of \cite{San2} was to solve the equivalence problem of  $\Delta GL(n;\mathbb R )$-structures,
finding the generators of rings of differential invariants of the structure (this is a general technique for
$G$-structures; see details in \cite{Jaime} y \cite{Valdes}). The construction of these differential invariants
leads to functorial connections attached to the $G$-structure which allows to construct differential invariants
in a natural way. About the rings of differential invariants of $\Delta GL(n;\mathbb R )$-structures we can
establish the following result:

\begin{theorem}
\emph{\cite{San2}} Let $r \in \mathbb N$, $r\geq 1$. The rings of
differential invariants of $r$-order the $\Delta GL(n;\mathbb R
)$-structures of a manifold $M$ are locally differentiably
generated over an open dense subset of $J^r(FM/\Delta GL(n;\mathbb
R ))$ by exactly $N_{2n,r}$ functions, where
\begin{equation}
N_{2n,r}=2n+
\binom{2n+r}{r}
 \frac{(3r-1)n^2-2(r+1)n}{r+1},
\label{eq:formula}
\end{equation}
If $r=0$ then $N_{2n,0}=0$; { i. e.}, the constant functions are the unique differential invariants of $0$-order.
\end{theorem}

In the case of $n=1$, the concept of the $\Delta GL(2; \mathbb R)$-structure over $M$ coincides with that of
$\mathbb R^*$-structure (see \cite{Valdes2}). Both types of $G$-structures define a $3$-web over a surface. By
the formula (\ref{eq:formula}) one has:
\[N_{2,r}=\frac{(r+1)(r-2)}{3}.\]

\noindent thus re-obtaining the result obtained by Vald\'{e}s in \cite{Valdes2}.

The explicit expression of a local basis of functions of the rings
of differential invariants of the $\Delta GL(n;\mathbb
R)$-structures is still an open problem. We hope to obtain
such invariants from the canonical connection of an almost
biparacomplex structure in a future work.

\subsection{Uniqueness of the functorial connection in the integrable case}

If $\mathrm{dim} M=2$, an $\alpha$-structure is always integrable
thus defining a web. The canonical connection is the Blaschke's
connection. In the paper \cite{MV} the authors have proved that
Blaschke's connection is the only functorial connection which can
be attached to two-dimensional three-webs. An open problem
consists on proving the following

\vspace{3mm}

{\bf Conjecture} In the biparacomplex case there exists only one functorial connection.

\subsection{Biparacomplex structures on Lie algebras and groups }

Let $\mathfrak{g}$ be a Lie algebra.  An   \emph{almost complex structure} $J$  on $\mathfrak{g}$ is a linear
endomorphism such that $J\circ J=-I$, where $I$ stands for the identity map. The structure is said {\em
integrable} or \emph{complex structure} if the corresponding Nijenhuis-type operator vanishes:
\begin{eqnarray*}
J[X,Y]=[JX,Y]+[X,JY]+J[JX,JY]
\end{eqnarray*}
for all $X,Y\in \mathfrak{g}$, where $[\, , \, ]$ denotes the Lie
bracket of the Lie algebra.

An {\em almost product structure}  $E$ on $\mathfrak{g}$ is a linear endomorphism such that $E\circ E=I$. The
structure is said {\em integrable} or \emph{product structure} if the corresponding Nijenhuis-type operator
vanishes:
\begin{eqnarray*}
  E[X,Y]=[EX,Y]+[X,EY]-E[EX,EY]
\end{eqnarray*}
This condition is equivalent to $\mathfrak{g}_{+}$ and $\mathfrak{g}_{-}$ being subalgebras, $\mathfrak{g}_{\pm
}$ being the eigenspace corresponding to the eigenvalue $\pm 1$ of the  product structure $E$.  Then
$(\mathfrak{g},\mathfrak{g}_{+},\mathfrak{g}_{-})$ is a double Lie algebra.

Observe that a Lie algebra $\mathfrak{g}$ is a real vector space $\mathbb{R}^{n}$ endowed with a Lie bracket
$[\, , \, ]$. As a real manifold, the Lie derivative of vector fields vanish, i.e., the Lie bracket of vector
fields vanish. Then, an almost complex (resp. product) structure is always integrable. But this is not the
situation for the Lie bracket of the Lie algebra.

One can recover some results of Andrada and Salamon \cite{AS}.
They consider a Lie algebra $\mathfrak{g}$ endowed with a pair $\{
J,E\} $ where $J$ is a complex structure on $\mathfrak{g}$, $E$ a
product structure on $\mathfrak{g}$ and $J\circ E=-E\circ J$.  Of
course, $P=J\circ E$ is also a product structure on
$\mathfrak{g}$, and $(\mathfrak{g},J,E,P)$ is a biparacomplex
manifold. They characterize double Lie algebras
$(\mathfrak{g},\mathfrak{g}_{+},\mathfrak{g}_{-})$ which are
associated to a complex product structure \cite[Prop. 2.5]{AS} and
prove that the complexification of a Lie algebra endowed with a
complex product structure has a hypercomplex structure \cite[Th.
3.3]{AS}. Explicit examples of complex product structures on
4-dimensional Lie algebras are given. All the results through the
paper of Andrada and Salamon are given only in the integrable
case. They prove \cite[Prop. 5.1]{AS} that a Lie algebra carrying
a complex product structure admits a unique torsion-free
connection parallelizing $J$ and $E$. This can be obtained as a
consequence of Theorems \ref{teor:cruceanu} and
\ref{teor:fn-torsion} of our paper.

Four-dimensional Lie algebras admitting a biparacomplex structure has been recently classified by Bla\v{z}i\'{c}
and Vukmirovi\'{c} \cite{BV} and by  Andrada,  Barberis,  Dotti and Ovando \cite{ABDO}. The authors of the first
paper use the name para-hypercomplex for such a structure.

Finally, we point out the work \cite{IT}, where the authors
consider a Lie group endowed with a biparacomplex structure
invariant respect to left translations. Then the Lie group is said
to admit a homogeneous complex product structure. They prove that
the Lie groups $SL(2m-1,\mathbb{R})$ and $SU(m,m-1)$ admit
homogeneous product structures.

The authors of the present paper think that the results obtained through the paper about almost complex product
structures on manifolds can be translated to the study of non-integrable complex product structures on Lie
algebras.

\subsection{Triple structures}
Almost biparacomplex manifolds are example of {\em triple
structures}, \emph{i.e.}, of manifolds endowed with  three
$(1,1)$-tensor fields $F$, $P$ and $J$ satisfying

 \[ F^2=\pm Id , \, \, P^2=\pm Id, \, \, J=P\circ F, \,
\, P\circ F\pm F\circ P=0.\]

In fact, one can define four different triple structures, namely

\begin{itemize}

\item {\em Almost biparacomplex structure}:  $F^{2}=Id,\,
P^{2}=Id, P\circ F+ F\circ P=0 $.

\item {\em Almost hyperproduct structure}: $F^{2}=Id,\, P^{2}=Id,
P\circ F- F\circ P=0 $.

\item {\em Almost bicomplex structure}: $F^{2}=-Id,\, P^{2}=-Id,
P\circ F- F\circ P=0 $.

\item {\em Almost hypercomplex structure}: $F^{2}=-Id,\,
P^{2}=-Id, P\circ F+ F\circ P=0 $.

\end{itemize}

Almost hyperproduct  and almost bicomplex structures do not admit
functorial connections (see \cite{EST} for a proof), whereas
almost biparacomplex and almost hypercomplex ones do admit: those
studied in this paper for almost biparacomplex structures; Obata
connection in the hypercomplex case, being the unique torsion-free
connection parallelizing the structure. As we have said in the
above subsection, a biparacomplex structure on a Lie algebra
$\mathfrak{g}$ defines a hypercomplex structure on its
complexification $\mathfrak{g}^{\mathbb{C}}$. Moreover, the Obata
connection on $\mathfrak{g}^{\mathbb{C}}$ is flat iff the
canonical connection on $\mathfrak{g}$ is flat (see \cite[Cor.
5.3]{AS}).

\subsection{Biparacomplex metric structures}

Let $(M,F,P)$ be an almost biparacomplex manifold. One of us
 has defined four different kinds of metrics adapted to
the biparacomplex structure

\begin{definition} {\em (see \cite{San})}.  Let $(M,F,P)$ be a biparacomplex
manifold, and let $g$ be a pseudo-Riemannian metric on $M$. Then, $(M,F,P,g)$ is said to be a $(\varepsilon
_{1},\varepsilon _{2})$ {\em pseudo-Riemannian almost biparacomplex manifold}, where $\varepsilon
_{1},\varepsilon _{2}\in \{+,-\}$ according to the following relations:
$$g(FX,FY)=\varepsilon _{1}g(X,Y);\:  g(PX,PY)=\varepsilon _{2}g(X,Y).$$
\end{definition}

Each one of the four possibilities of the signs determines the
sign of $g(JX,JY)$. In the  cases $(+,-), (-,+)$ the metric $g$ is
neutral of signature $(n,n)$ and in the case $(-,-)$ is neutral of
signature $(2n,2n)$. In \cite{Bl}, $(-,-)$ pseudo-Riemannian
almost biparacomplex manifold are called {\em paraquaternionic
Hermitian}. In that paper the Bla\v zi\'c studies the
paraquaternionic projective space, which is an example of this
structure.

\vspace{3mm}

A \emph{hypersymplectic} \cite{Hitchin}, \emph{hyper-Hermitian}
\cite{ITZ} or \emph{neutral hyperk\"{a}hler} \cite{Kamada}
manifold is a $4n$-dimensional biparacomplex  manifold endowed
with a neutral metric of signature $(2n,2n)$ such that
$g(JX,JY)=g(X,Y),\; g(FX,FY)=-g(X,Y)$. With the above notation, it
is a $(-,-)$-metric biparacomplex manifold. Then, it is K\"{a}hler
and Ricci flat. Moreover, hypersymplectic structures are used in
string theories.

 Recently, Andrada \cite{A} has classified hypersymplectic structures on
four-dimensional Lie algebras and Andrada and Dotti \cite{AD} have studied hypersymplectic structures on
$\mathbb{R}^{4n}$, showing significant examples.

\medskip

On the other hand, a connection with torsion attached to a
$(-,-)$-metric biparacomplex manifold have been considered in
\cite{ITZ}, where the authors define a \emph{hyperparaK\"{a}hler
with torsion} as a $(-,-)$-metric $g$ such that there exists a
linear connection $\nabla$ satisfying the following relations:
$$\nabla g=\nabla F=\nabla P=\nabla J=0; \; T(X,Y,Z)\doteqdot g(T(X,Y),Z)=-T(X,Y,Z)  $$
where $T$ denotes the torsion tensor of $\nabla$. The last
relation can be read saying that the torsion tensor of type (0,3)
is totally skew-symmetric. Moreover, they obtain a lot of
examples.

\medskip

Finally, relationships between almost biparacomplex structures and
almost bi-Lagrangian ones have been found by the authors. Let us
remember that an almost bi-Lagrangian structure on a symplectic
manifold $(M,\omega )$ is given by two transversal Lagrangian
foliations ${\cal D} _{1},{\cal D} _{2}$. Equivalently, it is an
almost para-K\"{a}hler structure on $M$ (see \cite{ES} or
\cite{EST2} for the details). Then one can prove (see
\cite{ESDebrecen}, \cite{EST2})

\begin{proposition}
Let $(M,\omega ,{\cal D} _{1},{\cal D} _{2})$ be an almost bi-Lagrangian manifold and let $(M,F,g)$
be its associated almost para-K\"{a}hler structure. For each Riemannian metric $G$ such that ${\cal D} _{1}$ and
${\cal D} _{2}$ are $G$-orthogonal, we define the almost complex structure $J$ associated with $G$ and $\omega $
(i.e., $\omega (X,Y)=G(JX,Y)$). Then:

(1) $(M,F,P=J\circ F)$ is an almost  biparacomplex manifold;

(2) $(M,J,g)$ is a Norden manifold;

(3) $(M,F,G)$ is a Riemannian almost product manifold;

(4) $(M,F,P,g)$ is a $(-,+)$ pseudo-Riemannian almost biparacomplex manifold;

(5) $(M,F,P,G)$ is a $(+,+)$ Riemannian almost  biparacomplex  manifold.
\end{proposition}

Such a metric always exists: if $H$ is any Riemannian metric on
$M$, then one can define a new Riemannian metric $G$ by
$G(X,Y)=H(X,Y)+H(FX,FY)$ obtaining that   $(M,F,G)$ is a
Riemannian almost product manifold, \emph{i.e.}, the two
distributions ${\cal D} _{1}$ and ${\cal D} _{2}$ are
$G$-orthogonal.
\vspace{3mm}

\noindent {\Large {\bf Acknowledgments}} \vspace{3mm}

F. E. is partially supported by   Project MTM2005-00173 (Spain) and R. S. is partially supported by MTM2005-05207 (Spain). The authors are also grateful to Professors V. Cruceanu, M. Fioravanti, P. M.
Gadea, J. Mu\~{n}oz Masqu\'{e}, and L. Ugarte for their comments
about some aspects of the topic. The second author wishes to thank
to his  thesis advisors, Prof. F. Etayo and Prof. J. Mu\~{n}oz
Masqu\'{e}, their valuable advice and infinity patience.

\end{document}